\documentclass{article}
\usepackage[utf8]{inputenc}
\usepackage{multicol}
\usepackage{amsmath}
\usepackage{amssymb}
\usepackage{amsthm}
\usepackage{hyperref}
\usepackage{enumerate}
\usepackage[letterpaper,margin=1.5in]{geometry}

\newtheorem{theorem}{\bf Theorem}[section]

\newtheorem{corollary}{\bf Corollary}[section]

\newtheorem{Proposition}[theorem]{Proposition}

\newcommand{\bomega}{\boldsymbol{\omega}}
\newcommand{\blambda}{\boldsymbol{\lambda}}
\newcommand{\derm}[1]{\mathrm{Der}^{#1}\,\Omega(M)}
\newcommand{\produ}[3]{\left\langle #1,#2;#3\right\rangle}
\newcommand{\produc}[2]{\left\langle #1;#2\right\rangle}
\def\nnabla{\nabla \hskip-2.45mm \nabla}

\begin{document}
\title{Supermanifolds, symplectic geometry and curvature}
\author{Rosal\'ia Hern\'andez-Amador$^*$, Jos\'e A. Vallejo$^\dagger$
and Juan Monterde$^\ddagger$ \\
{\normalsize $^*$Departamento de Matem\'aticas, Universidad de Sonora (M\'exico)}\\
{\normalsize $^\dagger$Facultad de Ciencias, Universidad Aut\'onoma de San Luis Potos\'i (M\'exico)}\\
{\normalsize $^\ddagger$Departament de Geometria i Topologia, 
Universitat de Val\`encia (Spain)}\\
{\footnotesize Emails: \texttt{guadalupehernandez@correoa.uson.mx,jvallejo@fc.uaslp.mx,juan.l.monterde@uv.es}}
}
\date{\today}
\maketitle

\begin{abstract}
We present a survey of some results and questions related to the 
notion of scalar curvature in the setting of symplectic supermanifolds.
\end{abstract}
\medskip
\begin{center}
\emph{Dedicado a Jaime Mu\~noz-Masqu\'e, maestro y amigo, en su 65 aniversario}
\end{center}

\section{Introduction}

Supermanifolds appeared in Mathematics as a way to unify the description of
bosons and fermions in Physics. Of course, there would be nothing special about
them if the resulting theory were just the juxtaposition of separate theorems, what is really interesting is the possibility of new phenomena arising from the interaction of both (the bosonic and the fermionic) worlds. From the point of view of Physics, the most prominent exponent is the phenomenon of supersymmetry, much questioned these days in view of the absence of experimental evidence coming from the LHC research, but from a purely mathematical point of view there is the exciting possibility of investigating geometric structures which can be understood only by looking at them through ``fermionic lenses''. 

Symplectic scalar curvature is one of these structures: if one starts out with a connection on a usual manifold, it is straightforward to define its associated curvature, but if a refinement such as Ricci or scalar curvature is desired (as in General Relativity), then a non-degenerate bilinear form (a second-order covariant tensor field) is required to take the relevant traces. Riemannian geometry enters the stage when that tensor field is taken symmetric, leading to a plethora of well-known results, but there is another possibility. A symplectic form could be used to make the successive contractions needed to pass from the curvature four-tensor to the scalar curvature, but it is readily discovered that the would-be symplectic scalar curvature  obtained this way vanishes due to the different symmetries involved (the Ricci tensor is symmetric and is contracted with the skew-symmetric symplectic form). Thus, it would seem that there is no room for a non-trivial Riemannian-symplectic geometry, an idea further supported from the observation that locally Riemannian and symplectic geometries are quite opposite to each other, as in the symplectic case there are no invariants because of the Darboux theorem.

However, things are different if we allow for supermanifolds. In this case, 
there are two variants of symplectic forms, even and odd ones, and it is 
remarkable that, while even symplectic forms lead to the same results as in the 
non graded setting, for odd symplectic manifolds it is possible, \emph{a priori},
to define a symplectic scalar curvature, because the symmetries involved in this
setting do not forbid its existence. However, the explicit construction of examples
is very difficult, and in this paper we try to explain why. The ultimate reason is
that the structure of odd symplectic manifolds is very restrictive. In particular,
they strongly depend on the existence on an isomorphism between the tangent bundle
$TM$ and the Batchelor bundle $E$ (that is, the vector bundle over $M$ such that
the supermanifold $(M,\mathcal{A})$ satisfies $\mathcal{A}\simeq \Gamma \Lambda E$).
When this isomorphism comes from a non-degenerate bilinear form on $TM$ with definite
symmetry (e.g, a Riemannian metric or a symplectic form), the symmetries of the
graded Ricci tensor lead to a trivial scalar curvature, as in the non-graded case.

While we will not deepen into the physical applications, neither of this odd
symplectic curvature nor supersymplectic forms in general
(for this, see \cite{BB081,BB082,CMQR08,RHA15}),
we will offer a detailed review of the mathematics involved in this construction
under quite general conditions, avoiding excessive 
technicalities with the aim of making this topic available to a wider audience.

\section{Preliminaries}\label{sec-prelim}

Let $M$ be a differential manifold, let $\mathcal{X}(M)$ denote the 
$\mathcal{C}^\infty(M)-$module of its vector fields, and let $\nabla$ be a
linear (Koszul) connection on it. The curvature of $\nabla$ is the operator
$\mathrm{Curv}:\mathcal{X}(M)\times \mathcal{X}(M)\to \mathrm{End}\mathcal{X}(M)$
such that
$$
\mathrm{Curv}(X,Y)=[\nabla_X,\nabla_Y]-\nabla_{[X,Y]}\,,
$$
where $[X,Y]$ is the Lie bracket of vector fields and $[\nabla_X,\nabla_Y]$ is
the commutator of endomorphisms.
Given a Riemannian metric on $M$ (that is, a
symmetric, positive-definite, covariant $2-$tensor field $g\in S^2_+(M)$), there
is a particular linear connection on $M$, the Levi-Civita connection, such that 
$\nabla g=0$. With the aid of the metric, two further contractions of the 
curvature can be defined, the first one leading to the Ricci covariant 
$2-$tensor
\begin{equation}\label{riccig}
\mathrm{Ric}(X,Y)=\mathrm{Tr}_g(Z\to \mathrm{Curv}(X,Z)Y)\,,
\end{equation}
and the second one to the Riemannian scalar curvature
\begin{equation}\label{scalarg}
S=\mathrm{Tr}(g^{-1}\mathrm{Ric})\,.
\end{equation}
Let us remark that the Ricci tensor \eqref{riccig} is symmetric, as it is $g$,
so the contraction in \eqref{scalarg} does not vanish \emph{a priori}.

Now suppose that we use a compatible symplectic form $\omega\in\Omega^2(M)$
(that is, such that $\nabla\omega =0$) to compute these contractions. Using a superindex to distinguish them from the previous ones, we obtain
\begin{equation}\label{ricciw}
\mathrm{Ric}^\omega(X,Y)=\mathrm{Tr}_\omega(Z\to \mathrm{Curv}(X,Z)Y)\,,
\end{equation}
and
\begin{equation}\label{scalarw}
S^\omega=\mathrm{Tr}(\omega^{-1}\mathrm{Ric}^\omega)\,.
\end{equation}
The symplectic Ricci tensor \eqref{ricciw} is again symmetric, but this time
the contraction in \eqref{scalarw} involves the skew-symmetric $\omega^{-1}$,
so we get $S^\omega =0$.

The study of symplectic manifolds $(M,\omega )$ endowed with a connection 
$\nabla$ such that $\nabla\omega =0$ can be carried on along lines similar to those of Riemannian geometry (see \cite{GRS98}). The resulting Fedosov manifolds
appeared first in the deformation quantization of Poisson manifolds (see 
\cite{FED94}). The fact that a basic local invariant such as the scalar curvature vanishes on any Fedosov manifold has led to a certain lack of interest in its use in Physics and Mathematics, aside from the mentioned r\^ole in
deformation quantization. However, if supermanifolds are considered a new
possibility appears. There are two classes of symplectic forms on a supermanifold and, as we see below, one of them has the symmetry properties
required to obtain a non-trivial contraction defining the symplectic scalar
curvature.

A supermanifold can be thought of as a non-commutative space of a special kind,
one in which the sheaf of commutative rings of $\mathcal{C}^\infty(M)$ 
functions has been replaced by a sheaf of $\mathbb{Z}_2-$graded supercommutative 
algebras, that is, to each open subset $U\subset M$ of a manifold, we assign an
algebra $\mathcal{A}(U)=\mathcal{A}_0(U)\oplus\mathcal{A}_1(U)$ with a product
such that 
$\mathcal{A}_i(U)\cdot\mathcal{A}_j(U)
\subset\mathcal{A}_{(i+j)\mathrm{mod}2}(U)$
and $a\cdot b=(-1)^{|a||b|}b\cdot a$, where $|a|$, $|b|$ denote the 
$\mathbb{Z}_2$ degree of the elements $a,b\in\mathcal{A}(U)$. An exposition of the basic facts about supermanifolds oriented to physical applications can be
found in \cite{SV09}. For completeness, let us give here the definition:
a real supermanifold is a ringed space $(M,\mathcal{A})$, where $\mathcal{A}$ is a sheaf of $\mathbb{Z}_2 -$graded commutative $\mathbb{R}-$algebras such that:
\begin{enumerate}[(a)]
\item If $\mathcal{N}$ denotes the sheaf of nilpotents of $\mathcal{A}$, then $\mathcal{A}/\mathcal{N}$ induces on $M$ the structure of a differential manifold.
\item The subsheaf $\mathcal{N}/\mathcal{N}^2$ is a locally free sheaf of modules, with $\mathcal{A}$ \emph{locally} isomorphic to the exterior sheaf $\bigwedge \left( \mathcal{N}/\mathcal{N}^2 \right)$.
\end{enumerate} 
The sheaf of differential forms on a manifold $M$, where $\Omega (U)=\bigoplus_{p\in \mathbb{Z}}\Omega^p (U)$, provide a good example.
The nilpotents in this case are all the $\alpha \in \Omega^p (M)$ with $p\geq 1$, so $\mathcal{A}/\mathcal{N}=\mathcal{C}^{\infty}(M)$ (the smooth functions on $M$). Moreover, $\mathcal{N}/\mathcal{N}^2=\Omega^1 (M)$, the space of $1-$forms, is locally generated by the differentials $\mathrm{d}x^1 ,...,\mathrm{d}x^m$ of the functions $x^i$ of a chart on $M$. Thus, as a model for a supermanifold we can think of a usual manifold $M$ endowed with ``superfunctions'', which are just differential forms and can be classified as even and odd by their degree. From now on, until otherwise explicitly stated, we will assume that our supermanifold
is $(M,\Omega (M))$, and sometimes we will refer to it as the Koszul or Cartan-Koszul 
supermanifold\footnote{This is not a great loss of generality in view of the existence
of the vector bundle isomorphism $TM\to E$, between $TM$ and the Batchelor bundle,
already mentioned in the Introduction (see \cite{Mon92}), 
so the changes needed to deal with
the most  general case are mainly notational.}.

The replacement of $\mathcal{C}^\infty(M)$ by
$\Omega (M)$ leads to the definition of other basic structures of differential
geometry. For instance, (super) vector fields on the supermanifold 
$(M,\Omega (M))$ are now the derivations $\mathrm{Der}\,\Omega (M)$ 
(such as the exterior differential $\mathrm{d}$, which has degree 
$|\mathrm{d}|=1$, 
the Lie derivative $\mathcal{L}_X$, which has degree $|\mathcal{L}_X|=0$, 
or the insertion 
$i_X$, which has degree $|i_X|=-1$). A straightforward corollary to a theorem of 
Fr\"ohlicher-Nijenhuis (see \cite{FN56}) states
that, given a linear connection $\nabla$ on $M$, the derivations of the form 
$\nabla_X$, $i_X$ generate the $\Omega (M)-$module $\derm{}$.

The (super) differential $1-$forms on 
$(M,\Omega (M))$ are defined as the duals $\mathrm{Der}^\ast \Omega (M)$, and
$k-$forms are defined by taking exterior products as usual, and noting that they
are \emph{bigraded} objects; if, for instance, 
$\boldsymbol{\omega}\in\Omega^2(M,\Omega (M))$ 
(that is the way of denoting the space
of $2-$superforms), its action on two supervector fields $D,D'\in\mathrm{Der}\,\Omega(M)$ will be denoted $\left\langle D,D';\boldsymbol{\omega}\right\rangle$, 
a notation well adapted to the fact that $\mathrm{Der}\,\Omega(M)$ is considered 
here as a left $\Omega(M)-$module and $\Omega^2(M,\Omega (M))$ as a right one. 
Other objects such as the graded exterior differential can be defined as in the 
classical setting, but taking into account the $\mathbb{Z}_2-$degree (for details in the spirit of this paper, see \cite{Val12}). Thus, if 
$\alpha\in\Omega^0 (M,\Omega(M))$, its graded differential $\boldsymbol{\mathrm{d}}$ 
is given by
$\left\langle D;\boldsymbol{\mathrm{d}}\alpha\right\rangle =D(\alpha)$, and if
$\boldsymbol{\beta} \in \Omega^1 (M,\Omega(M))$, we have a $2-$form 
$\boldsymbol{\mathrm{d}}\boldsymbol{\beta} \in \Omega^2 (M,\Omega(M))$ whose action is 
given by
$$
\left\langle D,D';\boldsymbol{\mathrm{d}}\boldsymbol{\beta}\right\rangle
=D(\left\langle D';\boldsymbol{\beta}\right\rangle)
-(-1)^{|D||D'|}D' (\left\langle D;\boldsymbol{\beta}\right\rangle)
-\left\langle [D,D'];\boldsymbol{\beta}\right\rangle\,,
$$
where $|D|$ denotes the degree of the derivation $D$.

\section{Symplectic supergeometry}

A supersymplectic form is a non-degenerate\footnote{In a technical sense that we will not describe here. See \cite{Kos77} for the details} graded $2-$form $\boldsymbol{\omega}\in\Omega^2(M,\Omega(M))$ such that $\boldsymbol{\mathrm{d}}\boldsymbol{\omega}=0$. Notice that there are two classes of supersymplectic forms: the even ones (for which $|\bomega|$ is even) act in such a way that, in terms of the induced $\mathbb{Z}_2-$degree,
$$
|\produ{D}{D'}{\bomega}|=|D|+|D'|
$$
and lead to symmetry properties similar to that of the non graded case,
but the odd symplectic forms (for which $|\bomega|$ is odd) satisfy
$$
|\produ{D}{D'}{\bomega}|=|D|+|D'|+1\,.
$$
As we will see below, these different properties translate into different symmetry
properties of the symplectic Ricci tensors.

By the aforementioned result of Fr\"olicher-Nijenhuis, given a linear connection
$\nabla$ on $M$, the study of the action of any $2-$superform $\bomega$ can be 
reduced to that of a matrix of the type
$$
\begin{pmatrix}
 \produ{\nabla_X}{\nabla_Y}{\bomega} & \produ{\nabla_X}{i_Y}{\bomega} \\
 \produ{i_X}{\nabla_Y}{\bomega} &  \produ{i_X}{i_Y}{\bomega}
\end{pmatrix}
$$
where $X,Y\in\mathcal{X}(M)$. 

In the case of an odd symplectic form $\bomega$, this structure can be made more 
explicit as follows. Starting from a vector bundle isomorphim $H:TM\to T^*M$, we
define an odd $1-$form $\blambda_H$, given by its action on basic derivations,
\begin{align*}
&\produc{\nabla_X}{\blambda_H}= H(X)\\
&\produc{i_X}{\blambda_H}= 0\,.
\end{align*}
(notice that this action is actually independent of $\nabla$). Next, we define 
$\bomega_H$ by $\bomega_H =\boldsymbol{\mathrm{d}}\blambda_H$. Thus,
the matrix of $\bomega_H$ now reads
\begin{align}\label{bomegah}
&\produ{\nabla_X}{\nabla_Y}{\bomega_H}= (\nabla_X H)Y-(\nabla_Y H)X\nonumber \\
& \produ{\nabla_X}{i_Y}{\bomega_H}= -H(X)(Y)\nonumber \\
& \produ{i_X}{\nabla_Y}{\bomega_H}= H(Y)(X) \\
& \produ{i_X}{i_Y}{\bomega_H}= 0\nonumber\,.
\end{align}
In a sense, these are all the odd symplectic superforms, 
according to the following result.
\begin{theorem}\cite{KM}
Let $\bomega$ be an odd symplectic form on $(M,\Omega(M))$,
then there exist a superdiffeomorphism $\phi:\Omega(M)\to
\Omega(M)$ and a fibre bundle isomorphism $H:TM\to T^*M$ such
that
$$\phi^\ast\bomega = \bomega_H\,.$$
\end{theorem}

In what follows, we will restrict our attention to \emph{odd} symplectic forms
of the type $\bomega_H$. Let us insist that the reason is that even symplectic
forms give rise to graded symmetric symplectic Ricci tensors (see \cite{RHA15}
for details), and further contraction with the graded skew-symmetric symplectic
form gives zero, thus leading to a trivial symplectic scalar supercurvature.

\section{Fedosov supermanifolds}
Now that we know the essentials about the structure of supersymplectic forms,
to begin the program sketched in Section \ref{sec-prelim} we
need some facts about superconnections $\nnabla$ on $(M,\Omega (M))$. 
In particular, we will need the analog of the Levi-Civit\`a theorem concerning
the existence of superconnections such that $\nabla\bomega =0$ for a 
supersymplectic form $\bomega$, and also their corresponding structure theorem.
We follow here the approach in \cite{MMV09}, although with some differences, the main
one being that we do not assume that $\nnabla$ is adapted to the splitting $H$
(also, see Theorem \ref{torsionless-nabla-compatible-odd-omega} below).

A superconnection on $(M,\Omega (M))$ is defined as in the non-graded case, as
an $\mathbb{R}-$bilinear mapping $\nnabla :\derm{}\times\derm{}\to\derm{}$,
whose action on $(D,D')$ is denoted $\nnabla_D D'$, with the usual properties of 
$\Omega (M)-$linearity in the first argument and Leibniz's rule in the second\footnote{In
particular, $\nnabla$ is not a tensor, hence the difference in notation.}:
$$
\nnabla_D(\alpha D') =D (\alpha )D' +(-1)^{|\alpha ||D|}\alpha\nnabla_D D'\,.
$$
The definition of torsion and curvature also mimics the non-graded case:
$$
\produ{D}{D'}{\mathrm{Tor}^{\nnabla}}=\nnabla_D D'-(-1)^{|D||D'|}\nnabla_{D'}D-[D,D']\,,
$$
and
$$
\produ{D,D'}{D''}{\mathrm{Curv}^{\nnabla}}=[\nnabla_D,\nnabla_{D'}]D''-
\nnabla_{[D,D']}D''\,,
$$
where $[D,D']=D\circ D'-(-1)^{|D||D'|}D'\circ D$,
$[\nnabla_D,\nnabla_{D'}]=\nnabla_D\nnabla_{D'}-(-1)^{|D||D'|}\nnabla_{D'}\nnabla_D$
are the graded commutators. As in the case of supersymplectic forms, 
we can describe a superconnection, once a linear connection $\nabla$ on $M$ is chosen,
by a set of tensor fields characterizing its action on basic derivations,
\begin{align*}
&\nnabla_{\nabla_X}\nabla_Y=\nabla_{\nabla_XY+K_0(X,Y)} + i_{L_0(X,Y)}\\[3pt]
&\nnabla_{\nabla_X}i_Y =\nabla_{K_1(X,Y)} + i_{\nabla_X Y+L_1(X,Y)}\\[3pt]
&\nnabla_{i_X}\nabla_Y=\nabla_{K_2(X,Y)} + i_{L_2(X,Y)}\\[3pt]
&\nnabla_{i_X}i_Y =\nabla_{K_3(X,Y)} + i_{L_3(X,Y)}\,,
\end{align*}
where $K_i,L_i :TM\otimes TM\to \Lambda T^*M\otimes TM$, for $i\in\{0,1,2,3\}$.
As a simplifying assumption, we will take a symmetric $\nnabla$. The relevant result
is the following.
\begin{theorem}\cite{MMV09}
Let $\nabla$ be a linear connection on $M$. A superconnection 
$\nnabla$ on $(M,\Omega (M))$ is symmetric if and only if
\begin{equation}\label{torsionless-connection}
\begin{array}{ll}
K_0(X,Y) = K_0(Y,X)-\mathrm{Tor}^\nabla(X,Y), & L_0(X,Y) =
L_0(Y,X)+ \mathrm{Curv}^\nabla(X,Y),\\
K_1(X,Y) = K_2(Y,X), & L_1(X,Y) = L_2(Y,X),\\
K_3(X,Y) = - K_3(Y,X), & L_3(X,Y) = - L_3(Y,X),
\end{array}
\end{equation}
for all $X,Y\in\mathcal{X}(M)$.
\end{theorem}
When the linear connection $\nabla$ on $M$ is symmetric, in the first equation of \eqref{torsionless-connection} we have,
$$
K_0(X,Y) = K_0(Y,X)\,,
$$
and this will be assumed in the sequel.

The next step is to study those superconnections $\nnabla$ which are compatible
with a given odd supersymplectic form $\bomega_H$, in the sense that 
$\nnabla\bomega_H =0$. This amounts to saying that
$$
D(\langle D_1,D_2;\bomega_H \rangle) = \langle \nnabla_D D_1,D_2;\bomega_H\rangle + (-1)^{|D||D_1|}\langle D_1,\nnabla_D D_2;\bomega_H\rangle\,, 
$$
for all $ D,D_1,D_2\in\derm{}$.
As a further simplifying assumption, we will take the linear
connection $\nabla$ compatible with the isomorphism $H:TM\to T^*M$, that is,
$\nabla H=0$ (so, \eqref{bomegah} also gets modified). Then, we get the following result (which corrects the one appearing in
\cite{MMV09}).
\begin{theorem}\label{torsionless-nabla-compatible-odd-omega}\cite{RHA15}
A symmetric superconnection, $\nnabla$, is compatible with the odd
symplectic form $\bomega_H$ if and only if
\begin{enumerate}[(a)]
\item\label{item1} $H(K_3(X,Y),Z)=-H(K_3(X,Z),Y)$
\item\label{item2-b} $H(K_2(X,Y),Z)=-H(Y,L_3(X,Z))$
\item $H(X,L_2(Y,Z))=H(Z,L_2(Y,X))$
\item\label{item4} $H(K_1(X,Y),Z)=H(K_1(X,Z),Y)$
\item $H(K_0(X,Y),Z)=-H(Y,L_1(X,Z))$
\item\label{item6} $H(X,L_0(Y,Z))= H(Z,L_0(Y,X))$,
\end{enumerate}
for all $X,Y,Z\in\mathcal{X}(M)$.
\end{theorem}
It is a straightforward generalization of the corresponding result in the non-graded 
setting, that superconnections compatible with a given supersymplectic form exist
and, moreover, they possess an affine structure (see \cite{RHA15} and, for a different approach \cite{Bla13}).
Also generalizing the non-graded case \cite{GRS98}, a Fedosov supermanifold is defined
as a supermanifold endowed with a supersymplectic form and a compatible symmetric 
superconnection, see \cite{GL04}. Combinig Theorems \ref{torsionless-connection} and
\ref{torsionless-nabla-compatible-odd-omega} with \eqref{bomegah}, we get the 
following.
Let $\bomega_H$ be an odd supersymplectic form on $(M,\Omega(M))$, with
$H:TM\to T^*M$ the associated bundle isomorphism. Let $\nabla$ be a compatible, 
symmetric, linear connection on $M$ (that is, $\nabla H=0$), so the action of
$\bomega_H$ on basic derivations reads
\begin{align}
& \produ{\nabla_X}{i_Y}{\bomega_H}= -H(X)(Y)\nonumber \\
& \produ{i_X}{\nabla_Y}{\bomega_H}= H(Y)(X)\label{bomegache} \\
& \produ{\nabla_X}{\nabla_Y}{\bomega_H}=0=\produ{i_X}{i_Y}{\bomega_H}\nonumber\,.
\end{align}
Finally, let $\nnabla$ be a superconnection on $(M,\Omega(M))$, symmetric and 
compatible with $\bomega_H$, characterized by the tensors $K_i,L_i$, $i\in\{0,1,2,3\}$.
Then, a pair $((M,\Omega (M)),\nnabla,\bomega_H )$ is a Fedosov 
supermanifold if and only if:
\begin{enumerate}[(a)] \setcounter{enumi}{6}
\item\label{itema} $K_0$ is symmetric, $L_0$ satisfies $L_0(X,Y)=L_0(Y,X)+\mathrm{Curv}^\nabla (X,Y)$,
and $K_3,L_3$ are skew-symmetric (from \eqref{torsionless-connection}).
\item\label{itemb} $K_1(X,Y)=K_2(Y,X)$ and $L_1(X,Y)=L_2(Y,X)$ (also from
\eqref{torsionless-connection}).
\item\label{itemc} The above items \eqref{item1} to \eqref{item6} hold.
\end{enumerate}
These conditions turn out to be very restrictive. From \eqref{item2-b},\eqref{itemb} and
\eqref{item4}, we get
$$
-H(X,L_3(Y,Z))=H(K_2(Y,X),Z)=H(K_1(X,Y),Z)=H(K_1(X,Z),Y)\,,
$$
and, because of the skew-symmetry of $L_3$ \eqref{itema}, this equals
$$
H(X,L_3(Z,Y))=-H(K_2(Z,X),Y)=-H(K_1(X,Z),Y)\,.
$$
Thus, $H(K_1(X,Z),Y)=-H(K_1(X,Z),Y)$,
which, in view of the fact that $H$ is an isomorphism, leads to 
$$
K_1=0=K_2
$$
and, \emph{a posteriori},
$$
L_3=0\,.
$$
An immediate consequence is the following.
\begin{corollary}\label{coro41}
A symmetric superconnection $\nnabla$, compatible with the odd symplectic form 
$\bomega_H$, acts as
\begin{align*}
&\nnabla_{\nabla_X}\nabla_Y=\nabla_{\nabla_XY+K_0(X,Y)} + i_{L_0(X,Y)}\\[3pt]
&\nnabla_{\nabla_X}i_Y = i_{\nabla_X Y+L_1(X,Y)}\\[3pt]
&\nnabla_{i_X}\nabla_Y=i_{L_1(Y,X)}\\[3pt]
&\nnabla_{i_X}i_Y =\nabla_{K_3(X,Y)}\,,
\end{align*}
for any $X,Y\in\mathcal{X}(M)$.
\end{corollary}
Notice that such a $\nnabla$ is determined just by four ordinary tensor fields
$K_0,K_3,L_0$, and $L_1$.

\section{Odd symplectic scalar curvature}
To study the simplest case,
we will start with an $n-$dimensional manifold $M$, an isomorphism 
$H:TM\to T^*M$ and a linear connection on $M$, $\nabla$, such that 
$\nabla H =0$.
We also consider the odd symplectic form $\bomega$ (actually $\bomega_H$,
but  we suppress subindices for simplicity) given by 
\eqref{bomegache} (denoting $H(X,Y)=H(X)(Y)$)
and a compatible superconnection $\nnabla$ as in corollary \ref{coro41}.
Due to the symmetry properties of $\mathrm{Curv}^{\nnabla}$, to characterize the action
of the symplectic curvature tensor
$$
\langle D_1,D_2,D_3,D_4 ; \mathbf R^{\bomega} \rangle 
   := \langle\, \langle D_1,D_2,D_3;\mathrm{Curv}^{\nnabla} \rangle \, , \, D_4 \,\, ; \, \bomega \rangle
$$
it suffices to study the following cases, which define corresponding $7$ tensor fields
$A_1,\ldots,A_5$, $B_1$, and $B_3$ (any other case gives a vanishing curvature) :
\begin{equation*}
 \begin{array}{rcl}
 \langle \nabla_X,\nabla_Y,\nabla_Z,\nabla_T \,;\, \mathbf{R}^{\bomega}\rangle 
      &=& H(T,B_1(X,Y,Z))  \\
 \langle \nabla_X,\nabla_Y,\nabla_Z,i_T \,;\,\mathbf{R}^{\bomega}\rangle       
      &=& -H(A_1(X,Y,Z),T) \\ 
      &=&  \langle \nabla_X,\nabla_Y,i_T,\nabla_Z \,;\,\mathbf{R}^{\bomega}\rangle\\
 \langle \nabla_X,i_Y,\nabla_Z,\nabla_T \,;\, \mathbf{R}^{\bomega}\rangle      
      &=& H(T,B_3(X,Y,Z)) \\ 
      &=&  - \langle i_Y,\nabla_X,\nabla_Z,\nabla_T \,;\,\mathbf{R}^{\bomega}\rangle \\
 \langle \nabla_X,\nabla_Y,i_Z,i_T;\mathbf{R}^{\bomega}\rangle 
      &=& -H(A_2(X,Y,Z),T) \\
 \langle \nabla_X,i_Y,\nabla_Z,i_T;\mathbf{R}^{\bomega}\rangle 
      &=& -H(A_3(X,Y,Z),T) \\
      &=& \langle \nabla_X,i_Y,i_T,\nabla_Z;\mathbf{R}^{\bomega}\rangle\\
      &=& -\langle i_Y,\nabla_X,\nabla_Z,i_T;\mathbf{R}^{\bomega}\rangle \\
      &=& - \langle i_Y,\nabla_X,i_T,\nabla_Z;\mathbf{R}^{\bomega}\rangle\\
 \langle \nabla_X,i_Y,i_Z,i_T;\mathbf{R}^{\bomega}\rangle 
      &=& -H(A_4(X,Y,Z),T) \\
      &=&   -\langle i_Y,\nabla_X,i_Z,i_T;\mathbf{R}^{\bomega}\rangle\\
 \langle i_X,i_Y,\nabla_Z,i_T;\mathbf{R}^{\bomega}\rangle 
      &=& -H(A_5(X,Y,Z),T) \\
      &=&   \langle i_X,i_Y,i_T,\nabla_Z;\mathbf{R}^{\bomega}\rangle\,.
 \end{array}
\end{equation*}
Of course, these new tensors can be explicitly computed from the $K_i,L_i$'s. For
instance, $A_2,A_3\in \Gamma (T^*M\otimes T^*M\otimes T^*M\otimes T^*M\otimes TM)$,
are given by 
\begin{align}
A_2(X,Y,Z)\cdot &= -K_3(\mathrm{Curv}^\nabla(X,Y)\cdot ,Z)\,\label{a2}\\
A_3(X,Y,Z)\cdot &= -K_3(Y,L_0(X,Z)\cdot )\,.\label{a3}
\end{align}
From these expression and items \eqref{item1}-\eqref{itemc} above, we get the following
\cite{RHA15}.
\begin{Proposition}
If $((M,\Omega (M)),\nnabla,\bomega )$ has the structure of a Fedosov supermanifold, then
\begin{enumerate}
\item $A_3(X,Y,Z)=A_3(Z,Y,X)-A_2(X,Z,Y)$.
\item $H(A_3(X,Y,Z),T)=H(A_3(Z,Y,X),T)-H(A_2(Z,X,T),Y)$.
\item $H(A_3(Y,Z,X),T)=-H(A_3(Y,T,X),Z)$,
\end{enumerate}
for any $X,Y,Z,T\in\mathcal{X}(M)$.
\end{Proposition}
If some additional symmetry properties of $H$ are added to these
conditions, we get those symmetries of the Ricci tensor mentioned in the introduction, leading to a trivial scalar curvature as we will see below.
\begin{corollary}\label{corosim}
If $H$ comes from a Riemannian metric or a symplectic form on $M$, then the graded
Ricci tensor satisfies
$$
\produ{\nabla_X}{i_Y}{\mathbf{Ric}^{\bomega}}=
-\produ{i_Y}{\nabla_X}{\mathbf{Ric}^{\bomega}}\,.
$$
\end{corollary}
Finally, we proceed to compute the symplectic scalar curvature from 
a graded Ricci tensor with this property. To this end, we take a basis
of homogeneous derivations $\{\nabla_{X_i},i_{X_i}\}$ (where $\{X_i\}$, for
$i\in\{1,\ldots,n\}$ is a local basis of vector fields on $M$). The odd supermatrix
locally representing $\bomega$ has the form
$$
\bomega =
\begin{pmatrix}
0 & -H(X_i,X_j) \\
H(X_j,X_i) & 0
\end{pmatrix} =
\begin{pmatrix}
0 & -H_{ij} \\
H^t_{ij} & 0
\end{pmatrix}\,.
$$
Thus, the graded morphism induced by $\bomega$, 
$\bomega^\flat :\derm{}\to\Omega^1(M,\Omega (M))$, has a supermatrix representative
$$
\bomega^\flat =
\begin{pmatrix}
0 & H_{ij} \\
-H^t_{ij} & 0
\end{pmatrix}\,.
$$
This supermatrix is invertible, and its superinverse is readily found to be
$$
(\bomega^\flat)^{-1}=
\begin{pmatrix}
0 & -(H^t_{ij})^{-1} \\
(H_{ij})^{-1} & 0
\end{pmatrix}\,.
$$
Now, the supermatrix associated to $\mathbf{Ric}^{\bomega}$ has the structure
$$
\mathbf{Ric}^{\bomega}=
\begin{pmatrix}
A & B \\
C & D
\end{pmatrix}\,,
$$
so
$$
(\mathbf{Ric}^{\bomega})^\flat =
\begin{pmatrix}
A^t & -(-1)^0 C^t \\
B^t & (-1)^0 D^t
\end{pmatrix}
=\begin{pmatrix}
A^t & -C^t \\
B^t & D^t
\end{pmatrix}\,.
$$
The scalar curvature is defined by the supertrace of $\mathbf{Ric}^{\bomega}$ 
with respect to $\bomega$; therefore, a straightforward computation shows that
\begin{align*}
\mathbf{Scal}^{\bomega} &= {\rm STr} \, \left( \left( \bomega^\flat\right)^{-1}\circ 
 \left(\mathbf{Ric}^{\bomega}\right)^\flat\right) \\
 &= -\mathrm{Tr}\left( C^t\, (H_{ij})^{-1} \right) 
 +\mathrm{Tr}\left( -B^t\, (H^t_{ij})^{-1} \right)\,.
\end{align*}
Now, if $H$ has a definite symmetry, from Corollary \ref{corosim} we get
$C=-B^t$ and consequently
$$
\mathbf{Scal}^{\bomega} =
-\mathrm{Tr}\left( C^t\, (H_{ij})^{-1} \right) 
 +\mathrm{Tr}\left( C\, (H^t_{ij})^{-1} \right)\,.
$$
But for any homogeneous invertible block $A$ we have
$$
(A^t)^{-1}=(-1)^{|A|}(A^{-1})^t
$$
(because, for homogeneous blocks, $(AB)^t=(-1)^{|A||B|}B^tA^t$), and also,
because of the invariance of the trace under transpositions,
$\mathrm{Tr}(A^t\,B)=\mathrm{Tr}(A\,B^t)$, so
$$
\mathbf{Scal}^{\bomega} = -\mathrm{Tr}(C\,(H^{-1}_{ij})^t)+ \mathrm{Tr}(C\,(H^{-1}_{ij})^t)=0\,.
$$
Thus, we deduce the following obstruction result (where we put back the subindex $H$ for
clarity).
\begin{theorem}
If $(M,H)$ is either a Riemannian or a symplectic manifold, then
$\mathbf{Scal}^{\bomega_H}=0$ on $(M,\Omega (M))$.
\end{theorem}

We believe that the preceding computations shed some light on the origin of 
the difficulties related to the construction of explicit examples of odd scalar
supercurvatures (letting aside the question of their geometric meaning).

Let us finish by mentioning two possible ways of avoiding this obstruction. 
Of course, one consists in taking a general $H:TM\to T^*M$, not symmetric 
nor skew-symmetric. The problem here is that such objects are not as natural
from the point of view of Physics as a metric or a symplectic form, and its
introduction should be carefully justified.
The other possibility involves the choice of a connection $\nabla$ such that
$\nabla H \neq 0$. This one is more interesting, as physically the choice of a 
connection is often part of the problem (for instance, in the Lagrangian version of 
Ashtekar's Canonical Gravity, connections are precisely the variables \cite{JS87}). 
However, the study of this case is much more difficult and will be treated somewhere
else \cite{RHA15}.


\thebibliography{C}

\bibitem{BB081} I. A. Batalin and K. Bering: \emph{Odd scalar curvature in field-antifield formalism}. J. Math. Phys. \textbf{49} (2008) 033515.
\bibitem{BB082} I. A. Batalin and K. Bering: \emph{Odd scalar curvature in anti-Poisson geometry}. Phys. Lett. B \textbf{663} 1-2 (2008) 132--135.
\bibitem{Bla13} P. A. Blaga: \emph{Symplectic connections on supermanifolds:
Existence and non-uniqueness}. Studia Univ. Babes-Bolyai Math. \textbf{58} 4 (2013) 477--483.

\bibitem{CMQR08} A. De Castro, I. Mart\'in, L. Quevedo and A. Restuccia:
\emph{Noncommutative associative superproduct for general supersymplectic forms}.
J. of High Ener. Phys. (2008) 0808:009.

 
\bibitem{FED94} B. V. Fedosov: \emph{A simple geometrical construction of deformation quantization}. J. of Diff. Geom. \textbf{40} 2 (1997) 213--238.

\bibitem{FN56} A. Frölicher, A. Nijenhuis: \emph{Theory of vector valued differential forms. Part I}. Indag. Math. \textbf{18} (1956) 338--360.

\bibitem{GRS98} I. Gelfand, V. Retakh and M. Shubin: \emph{Fedosov manifolds}.
Adv. in Math. \textbf{136} 1 (1998) 104--140.
\bibitem{GL04} B. Geyer and P. M. Lavrov: \emph{Fedosov supermanifolds: basic properties and the difference in even and odd cases}. Int. J. Mod. Phys.
\textbf{A 19} (2004) 3195.
\bibitem{RHA15} R. Hern\'andez-Amador. Ph. D. Thesis (Universidad de Sonora, M\'exico), in preparation (2015).

\bibitem{JS87} T. Jacobson and L. Smolin: \emph{The left-handed spin connection as a variable for canonical gravity}. Phys. Lett. \textbf{B196} (1987) 39--42.

\bibitem{KM} Y. Kosmann-Schwarzbach, J. Monterde, {\it Divergence
operators and odd Poisson brackets}. Ann. Inst. Fourier, {\bf 52},
419--456 (2002).
\bibitem{Kos77} B. Kostant: \emph{Graded manifolds, graded Lie theory, and prequantization}. In `Differential Geometrical Methods in Mathematical Physics'.
Lecture Notes in Mathematics \textbf{570} (1977) 177--306.
%
\bibitem{Mon92} J. Monterde: \emph{A Characterization of Graded Symplectic Structures}.  Diff. Geom. and its Appl. \textbf{2} (1992) 81--97.

\bibitem{MMV09} J. Monterde, J. Mu\~noz-Masqu\'e and J. A. Vallejo:
\emph{The structure of Fedosov supermanifolds}. J. of Geom. and Phys.
\textbf{59} (2009) 540--553.
\bibitem{SV09} G. Salgado and J. A. Vallejo: \emph{The Meaning of Time and Covariant Superderivatives in Supermechanics}. Adv. in Math. Phys. (2009)
Article ID 987524.
\bibitem{Val12} J. A. Vallejo: \emph{Symplectic connections and Fedosov's quantization on supermanifolds}. J. of Phys. Conf. Ser. \textbf{343} (2012) 012124.
\end{document}